\documentclass[12pt,a4paper]{article}

\usepackage{amssymb,amsmath,amsthm,a4wide,epsfig}
\begin{document}

\newcommand{\N}{\Bbb N}
\newcommand{\R}{\Bbb R}
\newcommand{\Z}{\Bbb Z}
\newcommand{\Q}{\Bbb Q}
\newcommand{\C}{\Bbb C}
\newcommand{\esp}{\vskip .3cm \noindent}
\def\CC#1{${\cal C}^{#1}$}
\def\h#1{\hat #1}
\def\t#1{\tilde #1}
\def\wt#1{\widetilde{#1}}
\def\wh#1{\widehat{#1}}

\newtheoremstyle{mypl}{10pt}{10pt}{\it}{}{\sc}{.}{ }{}
\theoremstyle{mypl}
\newtheorem{prop}{Proposition}
\newtheorem{lemma}[prop]{Lemma}
\newtheorem{cor}[prop]{Corollary}
\newtheorem{defi}[prop]{Definition}
\newtheorem*{thm}{Theorem}

\newtheoremstyle{myex}{10pt}{10pt}{\rm}{}{\sc}{.}{ }{}
\theoremstyle{myex}
\newtheorem*{ex}{Example}
\newtheorem*{exs}{Examples}
\newtheorem*{Qu}{Questions}
\newtheorem*{rem}{Remarks}
\newtheorem*{ack}{Acknowledgements}

\renewcommand\c{\mathcal{C}}
\renewcommand\d{\text{\sf diam}}

\title{Curves of constant diameter and inscribed polygons}
\author{Mathieu Baillif}
\date{\empty}
\maketitle
\abstract{We investigate the problem of finding simple closed curve $\Gamma$ in the Euclidean plane
that have the property $C_n(D)$ that on each $x$ in $\Gamma$, there is a regular $n$-gon with 
edges length $D$ inscribed in the 
curve at $x$, that is, with all its vertices lying on $\Gamma$, one vertex being $x$. If $n=2$, 
it is equivalent to having constant diameter, a property investigated by many great mathematicians, including
Euler, Hurwitz and Minkowski, for instance. We show using basic calculus and ideas
that dates back to Euler that there are $\c^{\infty}$ curves
satisfying $C_n(D)$, but that the circle is the only $\c^{2}$ regular curve satisfying $C_4(1/\sqrt{2})$ and
$C_2(1)$.
In an addendum, we improve this by showing that 
if a $\c^{2}$ regular curve satisfies $C_{2n}(D)$ and $C_2(R)$ for some $D,R>0$, then it is the circle.}

\section{Published Part}

{\em This section contains the text published in Elemente der Mathematik, 2009, vol. 64, no. 3. See the addendum for further results.}
\esp
\noindent
Let $\Gamma$ be a simple closed curve in the Euclidean plane. 
Say that a polygon $S$ is {\em inscribed in $\Gamma$ at $x$} if
all the vertices of $S$ lie on $\Gamma$ and one is $x$. 
A line segment is here considered as a $2$-gon. We say that $\Gamma$ has property
$C_n(D)$ (for some $n\ge 2$, $D>0$) if
\renewcommand\theequation{$C_n(D)$}
\begin{equation}
  \begin{array}{l}
  \forall x\in\Gamma \text{ there is a unique regular $n$-gon with edges length $D$}\\
  \text{inscribed in $\Gamma$  at $x$.}  \end{array}
\end{equation}
Notice that $C_2(D)$ is equivalent to the following since $\Gamma$ is simple and closed
($||\cdot ||$ is the Euclidean norm):
\renewcommand\theequation{$C(D)$}
\begin{equation}  \begin{array}{l}
  \forall x\in\Gamma \quad\exists ! y(x)\in\Gamma\text{ with }||x-y(x)||=D,\text{ and}   \\
  \text{if } z\not= y(x),\, ||x-z||<D.
  \end{array}
\end{equation} 

If one drops the unicity assumption, $C(D)$ is the property of having {\em constant diameter}, 
which is in fact
equivalent (for closed curves in the plane) to having 
{\em constant width} or {\em constant breadth} (for the definitions and the proof
of the equivalence, see \cite[chap. 25]{RademacherToeplitz}). 
It is a surprise to many (it was to me !) that curves of constant diameter different from 
the circle do exist.
The simplest examples (attributed to Reuleaux \cite{Reuleaux}, but implicit in a much earlier
paper by Euler \cite{EulerTriangularibus}) 
are pictured below. They are built with circle arcs whose
centers are marked with a black dot. Notice that these curves are not \CC{2}.
\begin{center}  
  \epsfig{figure=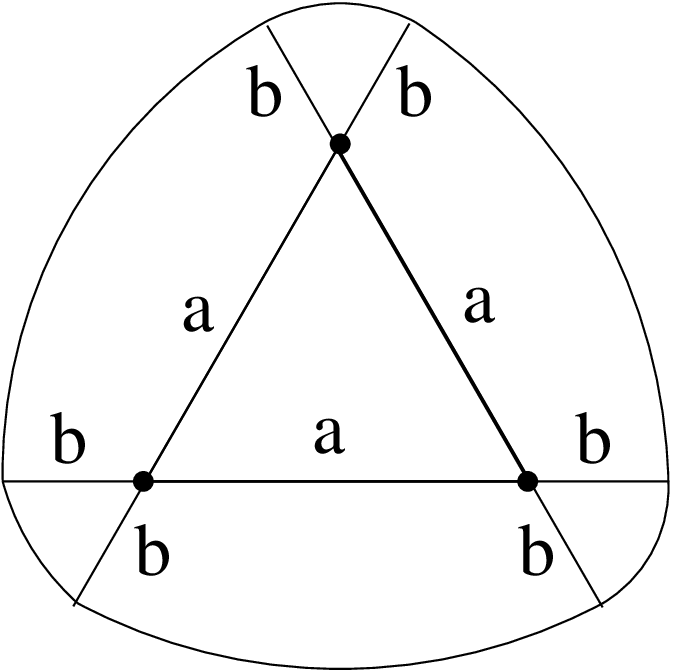,width=4cm}\hskip 1cm
  \epsfig{figure=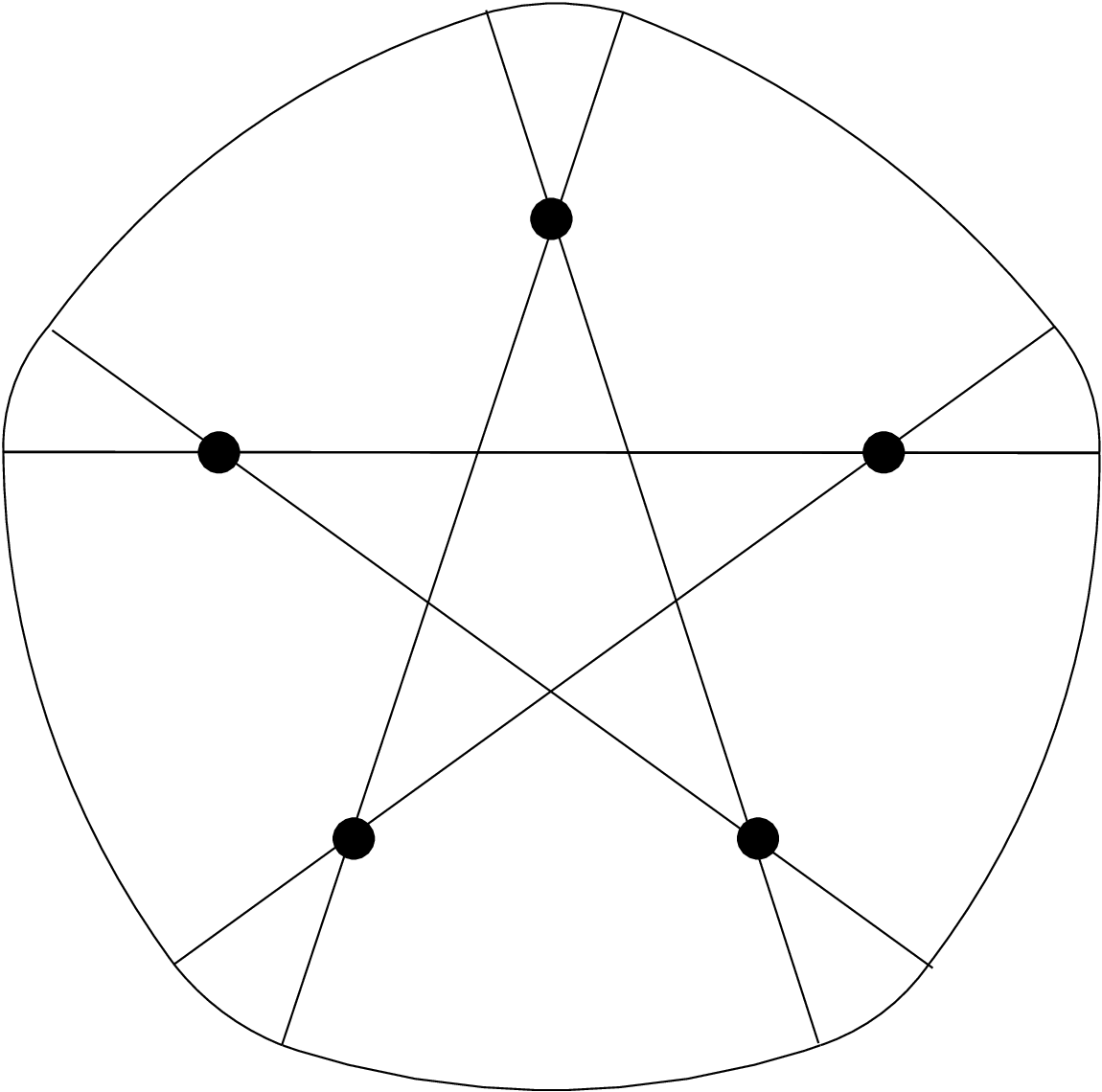,width=4cm}
\end{center}
The theory of curves of constant breadth has generated a considerable
literature, starting with Euler \cite{EulerTriangularibus} and including works by Hurwitz \cite{Hurwitz}, 
Minkowski \cite{Minkowski:konstanterBreite},
Blaschke (five articles in \cite{Blaschke:GesammelteWerke3}) and many others. 
These curves have many interesting properties, the most startling being perhaps that their
perimeter is $\pi$ times their diameter.
See \cite[chap. 25]{RademacherToeplitz} for an elementary account, or 
\cite{JordanFiedler} for a more thorough, though old, presentation.
For $n\ge 3$, to our knowledge the property $C_n(D)$ has not been investigated (but see 
\cite[p. 61]{JordanFiedler}, and
\cite{EgglestonTaylor} for a related problem).
\esp
In this note, we shall prove the theorem below,
using only
basic differential calculus.
We recall that a curve is {\em regular} if
it is \CC{1} with non vanishing derivative.

\begin{thm}     \begin{itemize}  
   \item[i) ] For all $n\ge 2$ and $D>0$, there are \CC{\infty} regular simple closed non       circular curves $\Gamma$ with property {\rm $C_n(D)$}.
  \item[ii) ] The circle of radius $D/2$ is the       only \CC{2} regular simple closed curve which satisfies both
      {\rm $C_4(\frac{D}{\sqrt{2}})$} and {\rm $C(D)$}.
  \end{itemize}
\end{thm}

\begin{rem}
  For i) $n=2$, such \CC{\infty} curves abound in the literature, see for instance \cite{Robertson} or
  \cite{Rabinowitz} for one given by a polynomial equation. 
  We however give a short self contained proof that yields
  simple explicit examples.
  Notice that ii) gives a definition of the circle involving only
  Euclidean distance between points {\em on} the circle, while most 
  usual definitions refer to points {\em off} the circle. 
\end{rem}

\proof[Proof of i) for $n=2$]
(Inspired by Euler \cite{EulerTriangularibus} and following G. Wanner's comments.)
To simplify we assume that $D=2$. 
We start we a stick of length $2$ which we place horizontally with center
$0$ at the origin. We then attach a needle at a point at position $-r$ from $0$ and we perform a small rotation
of angle $d\phi$ around this point. Then at any moment we change the position of the needle $r(\phi)$
and rotate the stick. As seen on the picture below, 
the center of the stick then describes a curve $(x(\phi),y(\phi))$ with
$$
  dx=-r(\phi)\sin\phi d\phi, \quad dy=r(\phi)\cos\phi d\phi,
$$
while the extremities of the stick follow the curve given by $$\gamma(\phi)=(x(\phi),y(\phi))+(\cos\phi,\sin\phi).$$

\begin{center}  \epsfig{figure=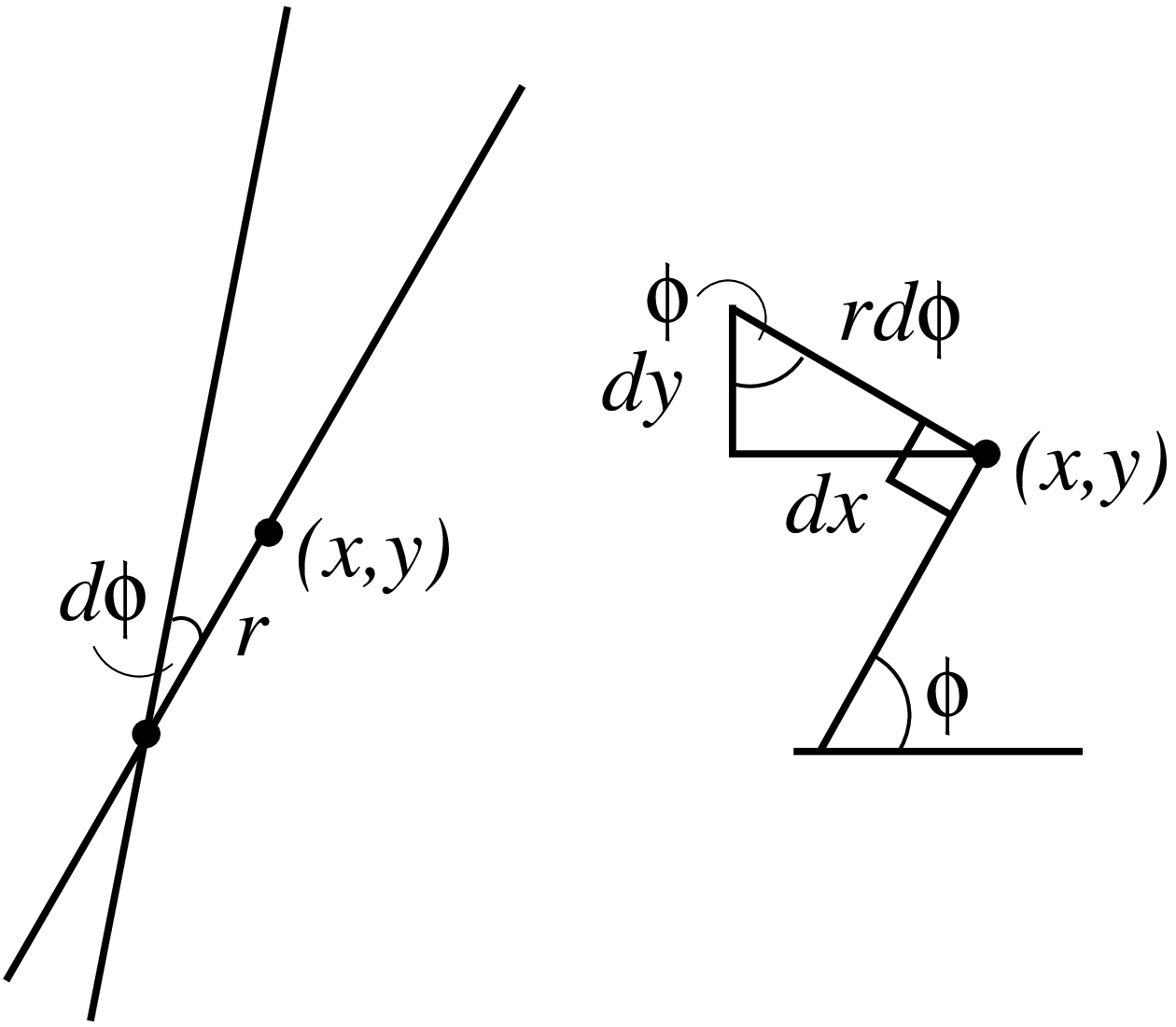,height=5cm}\hskip 1cm
\end{center}
When the angle $\pi$ is reached, the stick is again horizontal, and the position of the center is at the origin if $x(\pi)=y(\pi)=0$. The distance between $\gamma(\theta)$ and $\gamma(\phi)$ decreases (strictly if $|r(\theta)|<1$) as $\phi$ `goes away' from $\theta+\pi$, 
as seen on the figure below. (The circle of radius $2$ and center $\gamma(\theta)$ is dashed.) 
The curve obtained has thus constant diameter $2$ and satisfies {\rm $C(2)$} if $-1<r(\phi)<1$, strictly. Note that the proof works also if $r$ is only piecewise continuous. (A curve is piecewise continuous -- or \CC{1}, regular, etc -- if it is continuous -- or \CC{1}, regular, etc -- everywhere,  except possibly at a finite number of points.)
\begin{center}
  \epsfig{figure=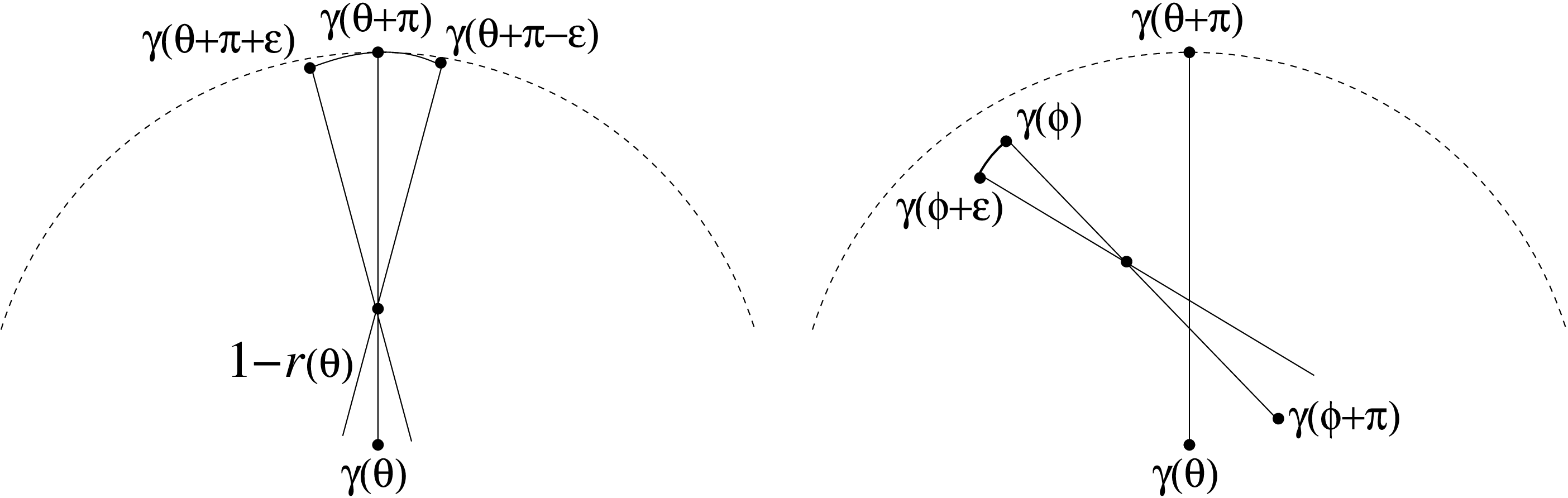,height=4cm}\hskip 1cm
\end{center}
To summarize, we proved:

\begin{prop}\label{prop1}
   Let $r:\R\to\R^2$ be piecewise continuous and       
   $$
     x(\theta)=-\int_0^\theta r(\phi)\sin(\phi)d\phi, 
     \quad y(\theta)=\int_0^\theta r(\phi)\cos(\phi)d\phi.
   $$
   Suppose that\\
   a) $r(\theta+\pi)=r(\theta)$, $|r(\theta)|\le\frac{D}{2}$,\\
   b) $(x(0),y(0))=(x(\pi),y(\pi))=(0,0)$.\\
   Then, the curve $\gamma(\theta)=(x(\theta),y(\theta))+\frac{D}{2}(\cos\theta,\sin\theta)$
   has constant diameter $D$, and
   has property {\rm $C(D)$} if $|r(\theta)|<\frac{D}{2}$.
\end{prop}
\endproof

Notice that given the path of the center $(x(\phi),y(\phi))$, the coordinate on the stick of the  `point of instant rotation' $r(\phi)$ is  given by the equation $x'(\phi)=-r(\phi)\sin(\phi)$.
The curves pictured in the first page of this note are obtained by taking $r$ piecewise constant. 
To obtain \CC{\infty} explicit examples, we can  
take $r(\theta)=a\cdot\sin((2k+1)\theta)$, with $k\ge 1$, $a<D/2$. Integrating, we get

\begin{equation*}
   (x(\theta),y(\theta))=\frac{a}{4}    \left(-\frac{\sin(2k\theta)}{k} + \frac{\sin(2(k+1)\theta)}{k+1} \, ,\,
    -\frac{\cos(2k\theta)}{k} -    \frac{\cos(2(k+1)\theta)}{k+1} \right).
\end{equation*}
By Proposition \ref{prop1}, the curve
$\gamma(\theta)=(x(\theta),y(\theta))+\frac{D}{2}(\cos\theta,\sin\theta)$ has property
$C(D)$.
One could also take any linear combination 
of $\cos((2k+1)\theta)$ and $\sin((2k+1)\theta)$ ($k\ge 1$), with small enough coefficients, for $r(\theta)$.
Below are pictured the curves of diameter $1$ (and $(x,y)$) for
$$    r(\theta)=\sin(3\theta)/3 + \cos(3\theta)/5,\, 
   \sin(5\theta)/2.01,\,
   \sin(3\theta)/10 + \cos(7\theta)/2.501.  
$$

\begin{center}
   \epsfig{figure=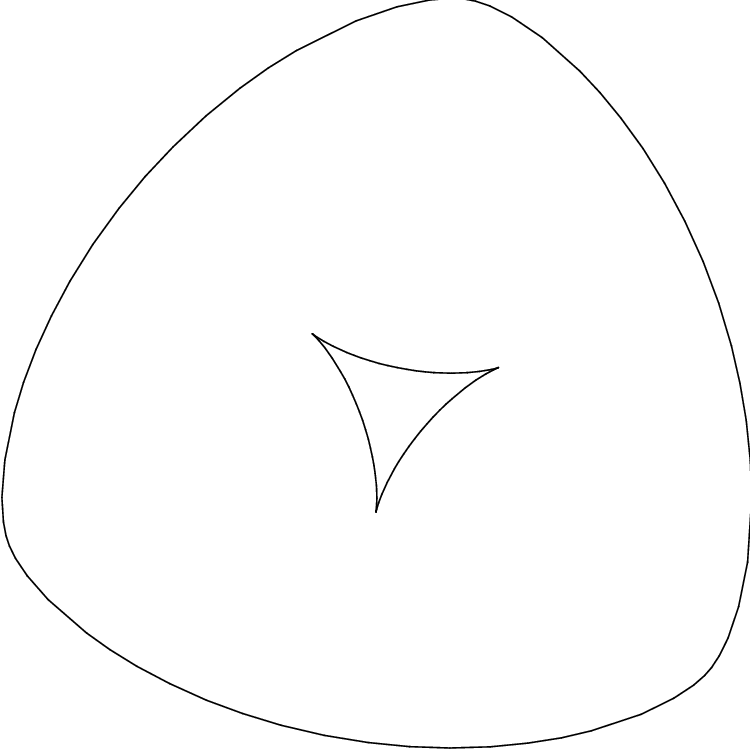,height=3cm}\hskip 1cm
   \epsfig{figure=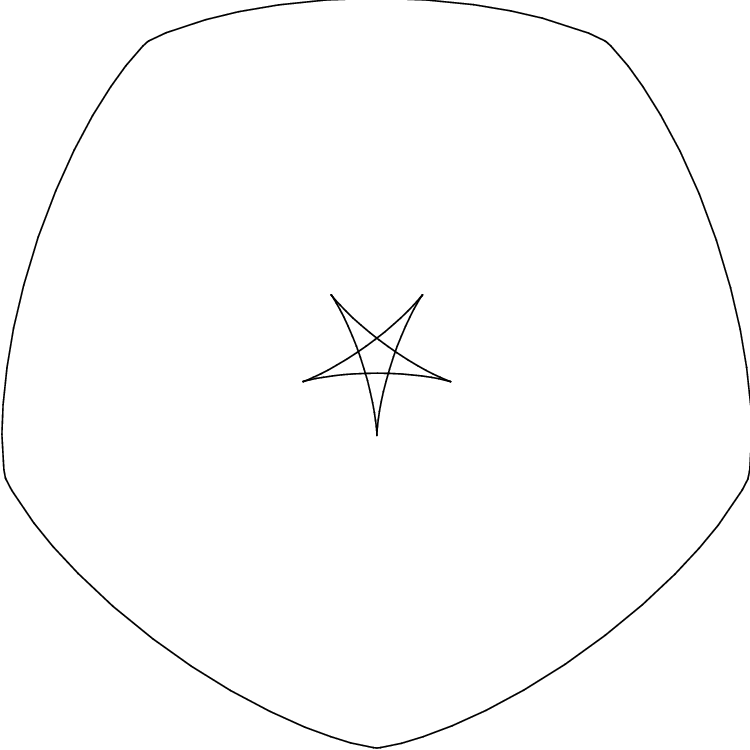,height=3cm}\hskip 1cm
   \epsfig{figure=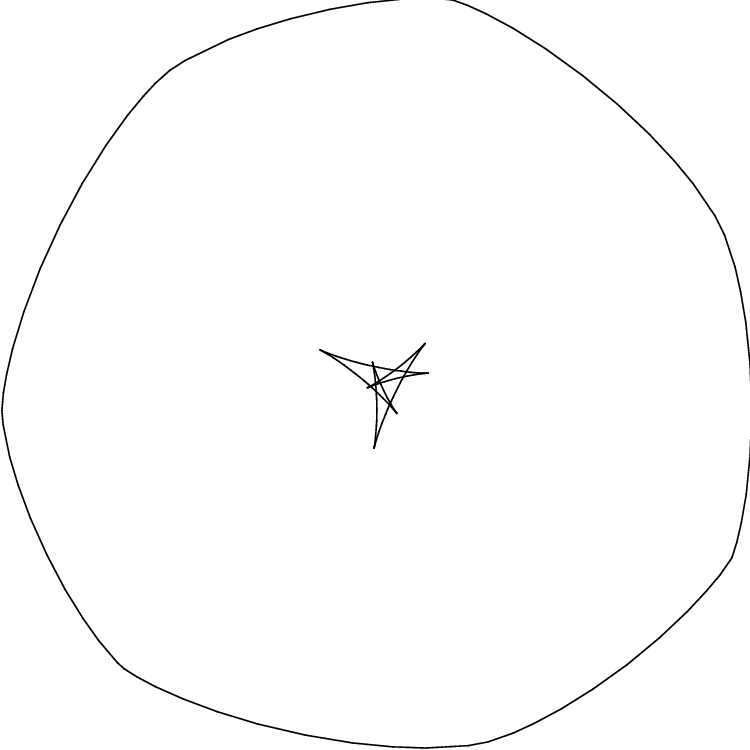,height=3cm}
\end{center}

\proof[Proof of i) for $n\ge 3$]
We use again G. Wanner's needle idea, but here we have more liberty. Take a regular $n$-gon which we rotate
with angle $\phi$ from $0$ to $2\pi/n$. At any moment we have {\em two} degrees of freedom to place our needle at a point inside the body of the $n$-gon, with coordinates $(\xi(\phi),\zeta(\phi))$ in a coordinate system moving with the $n$-gon. 
With the notations of the figure below,
the path $(x(\phi),y(\phi))$ of the center is determined by 
$dx=-r\sin(\phi + \psi)d\phi$, $dx=r\cos(\phi + \psi)d\phi$, and after expanding
$\sin(\phi + \psi)$ and $\cos(\phi + \psi)$ we get

$$  
  dx=-(\xi\sin(\phi)+\eta\cos(\phi))d\phi, \quad dy=(-\zeta\sin(\phi)+\xi\cos(\phi))d\phi.
$$

\begin{center}
   \epsfig{figure=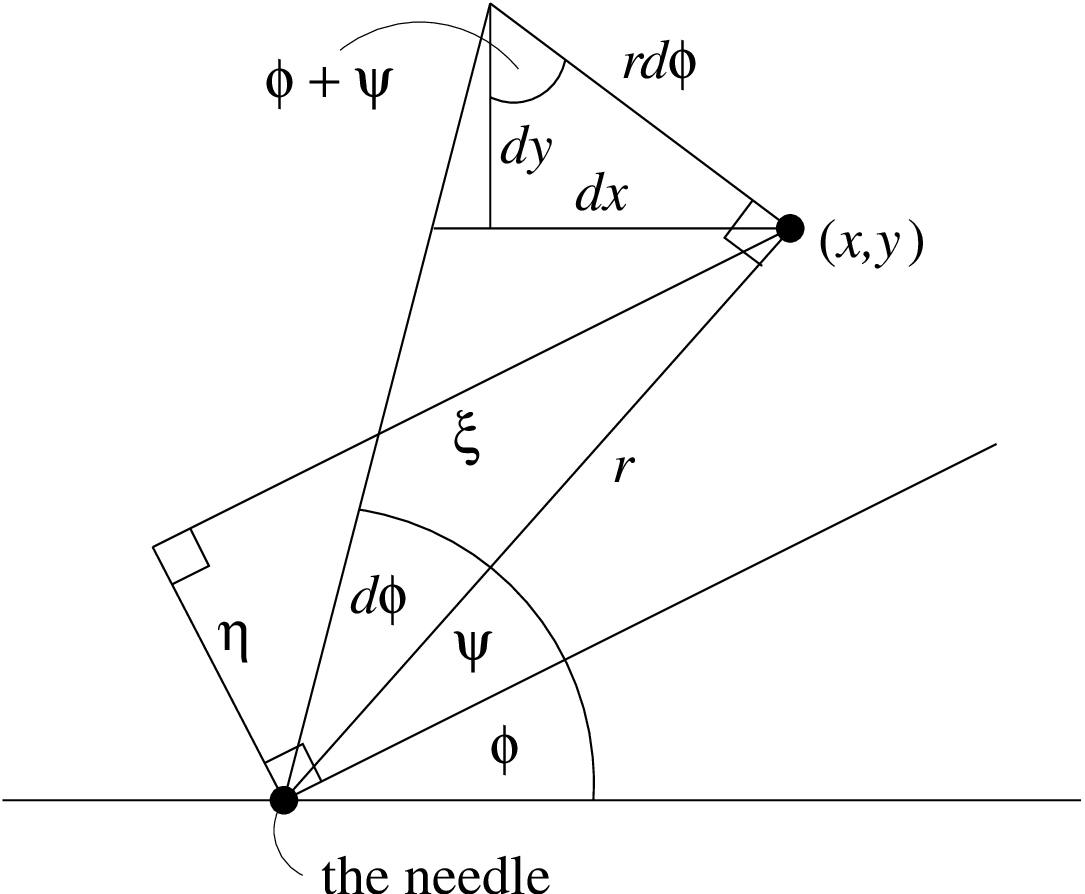, height=7cm}
\end{center}

Then, the condition for the curve described by the vertices of the $n$-gon to be closed is

$$   \int_0^{\frac{2\pi}{n}}(\xi(\phi)\sin(\phi)+\eta(\phi)\cos(\phi))d\phi=    \int_0^{\frac{2\pi}{n}}(-\eta(\phi)\sin(\phi)+\xi(\phi)\cos(\phi))d\phi=0.
$$
By construction there is of course an $n$-gon inscribed at each point,
and since there are exactly two points at distance $D$ from a given point, this $n$-gon is unique.

To obtain explicit examples it might be easier to {\em start} with $x(\phi),y(\phi)$ satisfying 
$x(\phi+\frac{2\pi}{n})=x(\phi)$, $y(\phi+\frac{2\pi}{n})=y(\phi)$ and let
$\gamma(\phi)=(x(\phi),y(\phi))+R(\cos(\phi),\sin(\phi))$, where $R$ is the radius of the $n$-gon with
edges length $D$. Then if $x,y$ and their derivatives are small enough 
so that we do not create `new' points at distance $D$ from a given point, we obtain a
curve satifying $C_n(D)$.
\endproof

\proof[Proof of ii)]
Let $\Gamma$ be a \CC{2} regular simple closed curve satisfying 
$C_4(D/\sqrt{2})$ and $C(D)$. We take $D=1$ for simplicity. 
Let $\gamma$ parametrise $\Gamma$ by arc length counterclockwise.\\
\parbox{0.55\textwidth}{If $x\in\Gamma$,   denote by $c(x)$ the unique point of $\Gamma$ for which $||x-c(x)||=1$. If $\gamma(t)=x$, 
  $\gamma'(t)$ must be   normal to $x-c(x)$, because the curve must be totally inside the ``dashed eye'' of the figure on the right, 
  and in fact, $c(x)$ is the point on the normal at $x$ at distance $1$ from $x$.}
\hfill
\parbox{0.4\textwidth}{\epsfig{figure=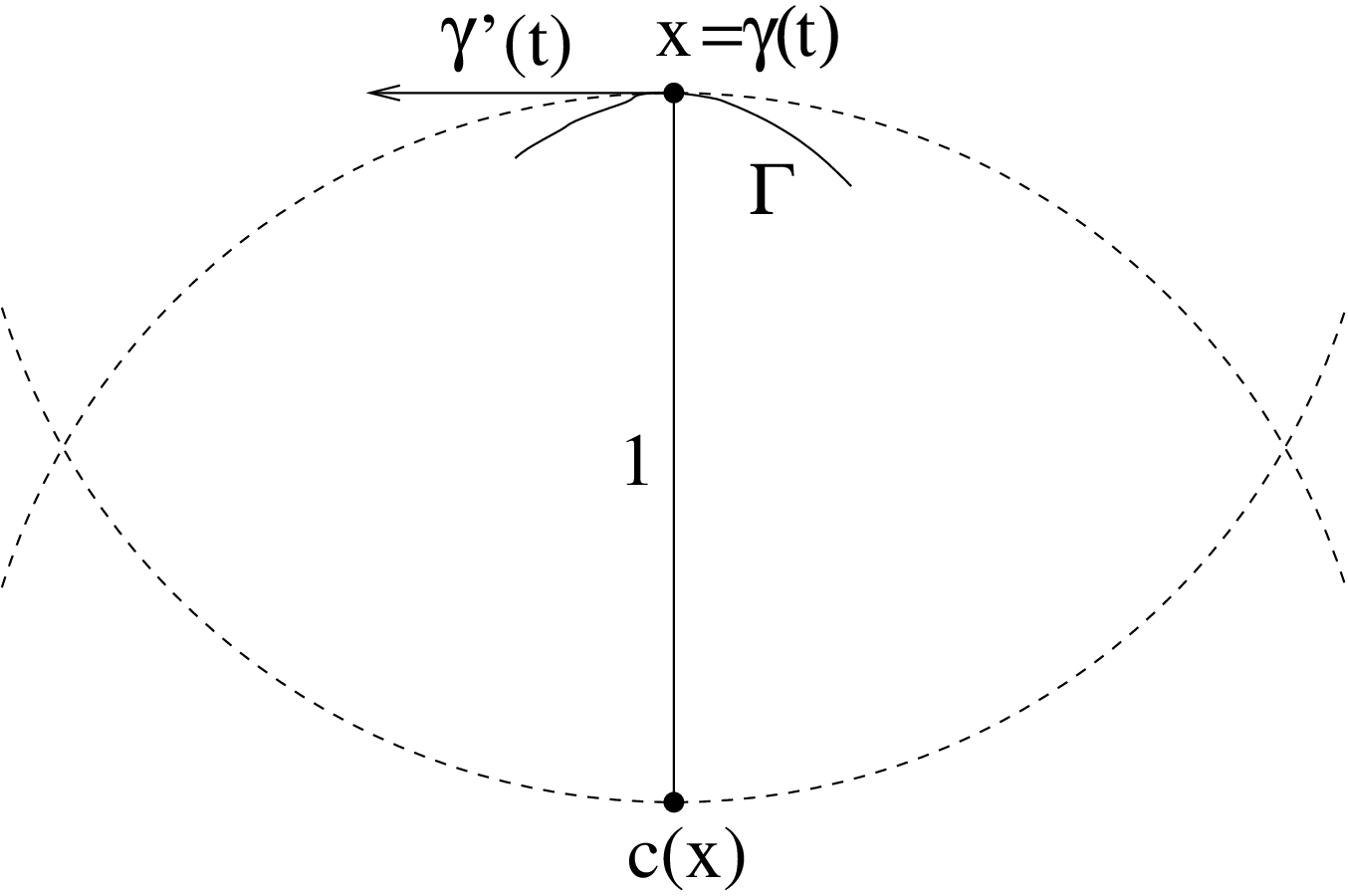,height=3.8cm}}
\\ 
Since the unit normal vector at $\gamma(t)$ is $n(t)=(-\gamma_2'(t),\gamma_1'(t))$, $c(\gamma(t))=\gamma(t)+n(t)$ is differentiable in $t$.
(Notice that we cannot use the implicit function theorem here, because the derivative of $||x-z||$ vanishes
precisely at $z=c(x)$.)
Moreover, the angle between the oriented line segment $[\gamma(t),c(\gamma(t))]$ 
and the horizontal axis strictly increases. 
Given $\theta\in[0,2\pi]$,
there is thus a unique (oriented) line segment $[x,c(x)]$ (the stick) making an angle $\theta$ with the horizontal axis,
we define $G(\theta)=(x(\theta),y(\theta))$ to be (the coordinates of) its middle point.
Since $c(x)$ is differentiable with respect to $x$, 
$G(\theta)$ is differentiable with respect to 
$\theta$.
Then, $\wt{\gamma}(\theta)=G(\theta)+\frac{1}{2}(\cos\theta,\sin\theta)$ is a parametrisation of $\Gamma$
(see the figure below, on the left).
By definition, $c(\wt{\gamma}(\theta))=\wt{\gamma}(\theta+\pi)$, and 
since the tangent of $\Gamma$ at $x$ is normal to $x-c(x)$, we have

\begin{equation*}
  \bigl<\wt{\gamma}'(\theta)\,|\,(\cos\theta,\sin\theta)\bigr> 
  =\bigl< G'(\theta)\,|\,(\cos\theta,\sin\theta)\bigr>=0,
\end{equation*}
which implies
$G'(\theta)=r(\theta)(-\sin\theta,\cos\theta)$ for some continuous $r$. 
(Recall that $r$ gives the coordinate on the stick of the point of instant rotation.) 
\\
Now, since $\Gamma$ has property $C_4(1/\sqrt{2})$, there is a unique square $S(\theta)$ 
with edges length $\frac{1}{\sqrt{2}}$ inscribed in $\Gamma$ at $\wt{\gamma}(\theta)$.
Since $\wt{\gamma}(\theta+\pi)$ is the unique point of $\Gamma$ at distance $1$ from 
$\wt{\gamma}(\theta)$, 
$\wt{\gamma}(\theta+\pi)$ is the vertex of $S(\theta)$ diagonal to $\wt{\gamma}(\theta)$ ($S(\theta)$ has
diagonal $1$). 
Thus, $G(\theta)$ is also the center of $S(\theta)$, which implies $G(\theta+\pi/2)=G(\theta)$ 
(see figure below, on the right).

\begin{center}
  \epsfig{figure=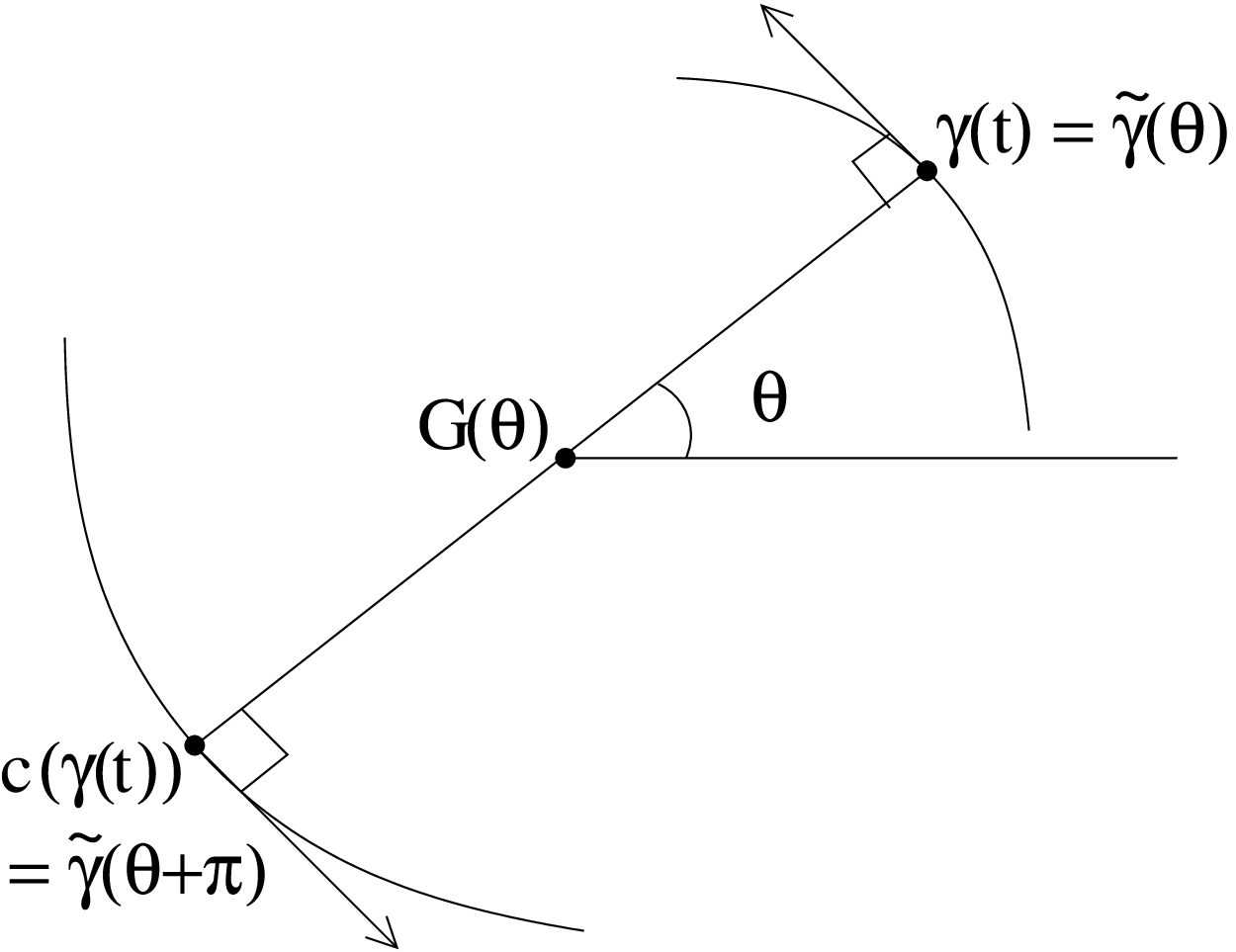,height=4cm}\hskip 1cm
  \epsfig{figure=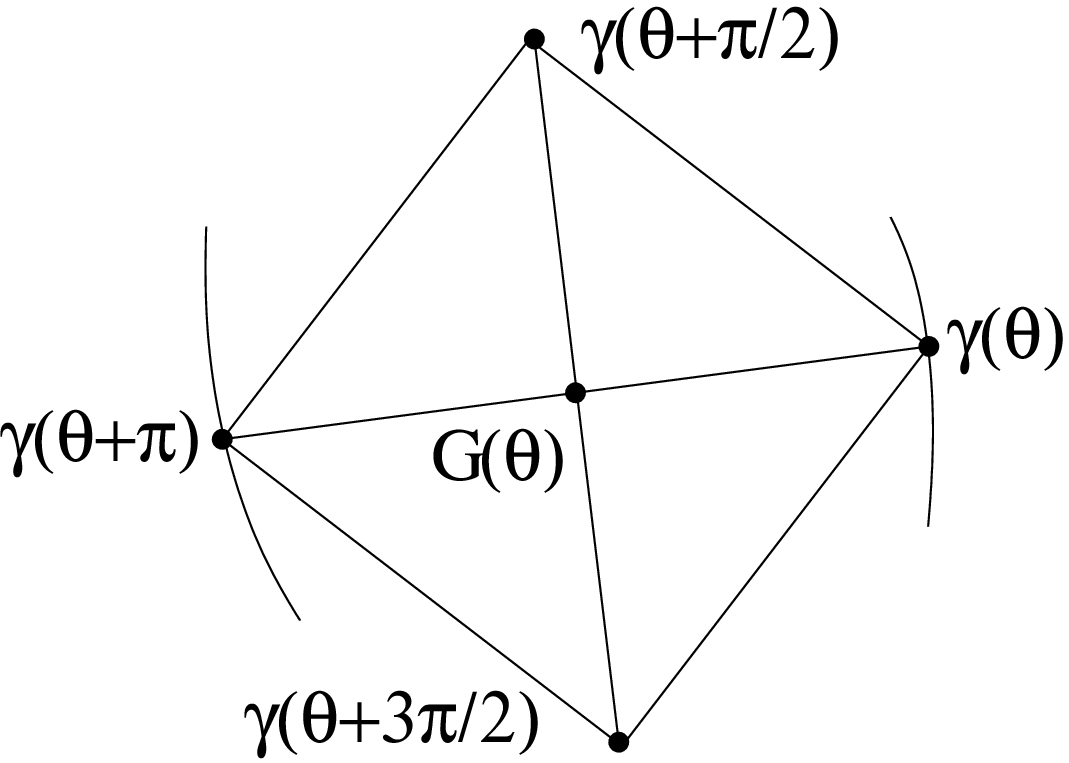,height=4cm}
\end{center}
Therefore, $G'(\theta+\pi/2)=G'(\theta)$, i.e.
$r(\theta+\pi /2)(\cos\theta,\sin\theta)=r(\theta)(-\sin\theta,\cos\theta)$,
so $r(\theta)=0$, $G(\theta)$ is constant, and hence $\Gamma$ is a circle.
\endproof

The end of the proof can also be seen as follows. 
Since we rotate a square, the needle 
must stay on the first diagonal of the square on the first half of the trajectory in order to have $C(1)$.
But on the second half it must stay on the second diagonal, and thus it must stay in the center
and the curve is a circle. 

This proof can be easily generalized to show that
if a \CC{2} regular closed simple curve has both properties 
$C_{2n}(D)$ and $C(R)$, with $R$ 
two times the radius of the regular $2n$-gon with edges length $D$, then it is the circle.
This however leaves open the following:

\begin{Qu}
   Are there $D,\wt{D}>0$ and 
   curves other than the circle which have property $C(D)$ and $C_n(\wt{D})$ for some $n\ge 3$~?
\end{Qu}

We finish with the following proposition, which gives
a motivation for the ``$!$'' in the definition of $C(D)$: 

\begin{prop}\label{prop2}
  Let $\Gamma$ be a continuous closed piecewise regular curve satisfying {\rm $C(D)$}. Then,
  $\Gamma$ is regular.
\end{prop}
The proof is an exercise for the interested reader.
({\em Hint: } Show that $c(x)$ is continuous in $x$, then that $\Gamma$ cannot have corners.) 
Notice that
taking $b=0$ in the first figure of this paper yields a curve of constant diameter with corners (the well
known Reuleaux triangle).

\begin{ack}
  I wish to thank G. Wanner, who gave so many suggestions for improvements that he should 
  be listed as a co-author,
  and B. Dudez, librarian at Geneva's math department,
  for his bibliographical help. I thank also D. Cimasoni.
\end{ack}

\section{Addendum: Improving Theorem 1 ii)}

The following is an easy improvement of Theorem 1 ii) and answers partially the questions above. 
It should have been guessed before submitting the paper, but unfortunately was not. 

\begin{prop}\label{propAdd}
   Let $\Gamma$ be a continuous simple closed curve satisfying 
   $C_{2n}(D)$ and $C(R)$. 
   Then $R$ is
   twice the radius of the regular $2n$-gon with edges length $D$.
\end{prop}

\proof
   Given $x\in\Gamma$, we denote the $2n$-gon inscribed at $x$ by $G(x)$.
   First, the vertices of $G(x)$ depend continuously on $x$. Indeed, given a sequence $x^1_m\in\Gamma$ ($m\in\N$)
   converging to $x^1\in\Gamma$, write $x^k_m$ for the sequences of the vertices of $G(x_m)$ ($k=1,\dots,2n-1$).
   Since $\Gamma$ is a closed curve, by compactness there is a subsequence $x^2_{i(k)}$ converging to some $x^2$, for
   some strictly increasing $i:\N\to\N$. Proceeding 
   by induction and taking subsequences at each step, we obtain a strictly increasing $j:\N\to\N$ with $x^k_{j(m)}$ converging
   to $x^k$, for each $k=1,\dots,2n-1$. By continuity, the $x^k$ are the vertices of a regular $n$-gon inscribed at $x^1$.
   Since this $n$-gon is unique, any converging subsequence of the $x^k_m$ have to converge to $x^k$, proving that the sequences
   $x^k_m$ themselves converge to it (using the compactness of $\Gamma$ again). This proves the continuity of these vertices.
   In a $2n$ gon, each vertex $v$ has a unique farthest vertex $o(v)$, so given $x\in\Gamma$, we denote also by $o(x)$ the 
   point of $\Gamma\cap G(x)$ farthest to $x$. Write $\theta_o(x)$ for the angle that the oriented line segment $[x,o(x)]$ makes with the
   horizontal axis, which varies continuously with $x$.

   The above arguments for $n=2$ show that the `antipodal point' $c(x)$ of $x$ depends continuously on $x$ as well. Moreover, as seen in the proof of 
   Theorem 1 ii), the angle of the oriented line segment 
   $[x,c(x)]$ increases strictly as $x$ follows $\Gamma$ counterclockwise, and depends continuously on $x$. We denote this angle by $\theta_c(x)$.
   
   Take $z\in\Gamma$ such that $\theta_c(z)=0$.      
   Notice that the length of the line segment $[x,o(x)]$ is twice the 
   radius of the $2n$-gon, so if
   $R$ is not equal to it,
   $\theta_o(x)\not= 0$.
   Consider the two arcs joining $z$ to $c(z)$. One of them must contain some $y$ with $\theta_o(y)=0$.
   Seeing this arc as an interval, we thus have (see the figure below, left)
   $$
      z \le y < o(y) \le c(z).
   $$
   But since $\theta_c(z)=\theta_o(y) = 0$, $\theta_c(c(z))=\theta_o(o(y)) = \pi$ and $\theta_c$ is strictly increasing, there is one
   point $x$ for which $\theta_c(x) = \theta_o(x)$ (see figure below, right, $\theta_c$ is in plain line and $\theta_o$ in dashed line).
   But then $c(x) = o(x)$, which shows that $R$ is equal to the length of the line segment $[x,o(x)]$, that is, to twice the radius of the $n$-gon.

\begin{center}
  \epsfig{figure=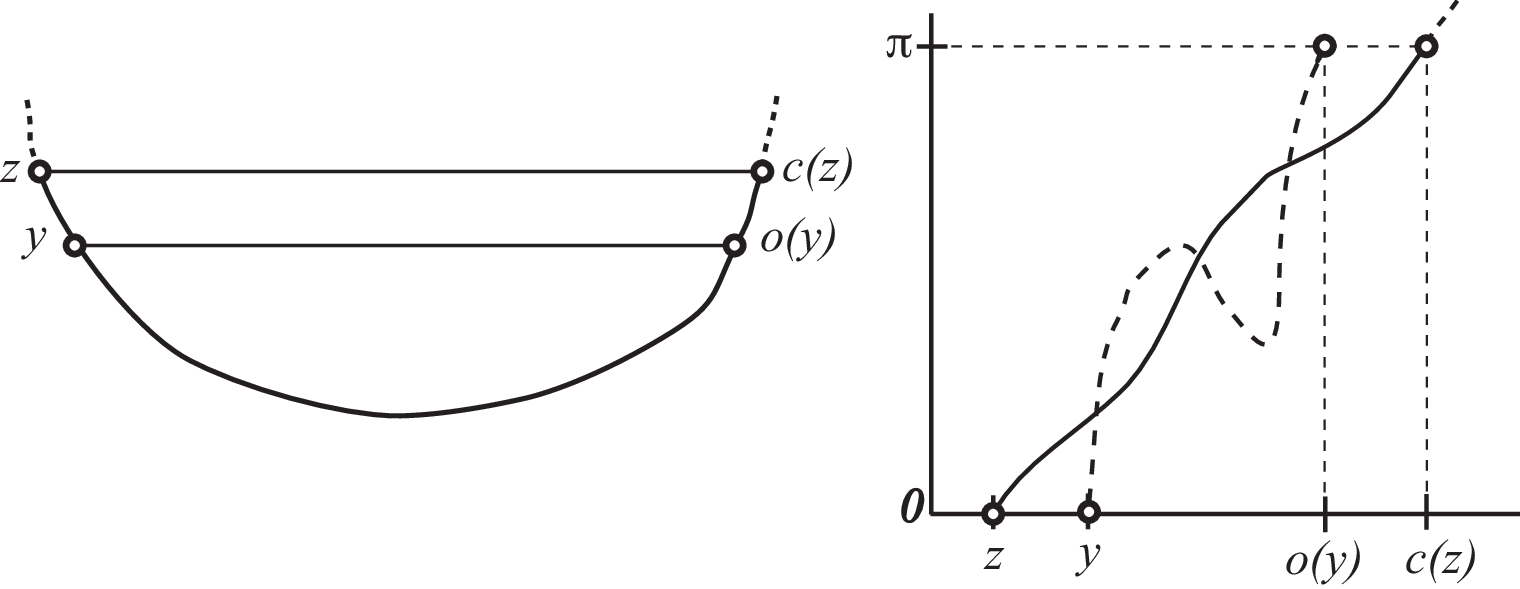,height=5cm}
\end{center}
  This finishes the proof.
\endproof

As an immediate corollary we get:

\begin{cor}
   If a $\c^{2}$ regular curve satisfies $C(R)$ and $C_{2n}(D)$ for some $R,D>0$, it is the circle.
\end{cor}
\proof
  By Proposition \ref{propAdd} and Theorem 1 ii).
\endproof

\end{document}